\documentclass{article}

\usepackage{amsmath}
\usepackage{amssymb}
\usepackage{amsthm}

\newtheorem{theorem}{Theorem}[section]
\newtheorem{lemma}[theorem]{Lemma}

\newcommand{\MC}{\mathcal{MC}}
\newcommand{\SR}{\mathcal{SR}}
\newcommand{\TA}{\mathcal{TA}}

\DeclareMathOperator\conv{conv}
\DeclareMathOperator\diag{diag}
\DeclareMathOperator\rk{rk}

\begin{document}

\title{Trigonometric approximation of the Max-Cut polytope is star-like}
\author{Romain Ageron \thanks{%
Paris-Saclay University, CentraleSup\'elec, 91190 Gif-sur-Yvette, France
({\tt romain.ageron@student-cs.fr}).}}

\maketitle

\begin{abstract}
The Max-Cut polytope appears in the formulation of many difficult combinatorial optimization problems. These problems 
can also be formulated as optimization problems over the so-called \textit{trigonometric approximation} which possesses 
an algorithmically accessible description but is not convex. Hirschfeld conjectured that this trigonometric approximation is 
star-like. In this article, we provide a proof of this conjecture.
\end{abstract}

Keywords: Max-Cut polytope, Trigonometric Approximation

MSC 2020: 90C20, 90C27

\section{Introduction}

A common problem in combinatorial optimization is the maximization of a quadratic form over \(\{-1,1\}^n\)
\begin{equation}
	\max_{x \in \{-1,1\}^n} x^T A x ~~= \hspace{-2ex} \max_{
		\scriptsize
		\setlength{\tabcolsep}{0pt}
		\begin{array}{c}
			X = x x^T \\ 
			x \in \{-1,1\}^n
		\end{array}
	} \langle A, X \rangle
	\label{eqn:problem}
\end{equation}
where \(\langle ., . \rangle\) denotes the usual scalar product on real symmetric matrices of size \(n\).

The decision problem associated to this optimization problem is NP-complete. Indeed the Max-Cut problem, one of 
Karp's 21 NP-complete problems, can be reduced in polynomial time to the maximization of a quadratic form over 
\(\{-1,1\}^n\) \cite{Hirschfeld}. The reformulation in the form of (\ref{eqn:problem}) of several common hard 
combinatorial optimization problems such as vertex cover, knapsack, traveling salesman, etc, can be found in 
\cite{Lucas_2014}.

Consider the set
\[
	\SR = \{ X \succeq 0 ~|~ \diag X = 1 \}
\]
in the space of real symmetric \(n \times n\) matrices, where \(X \succeq 0\) means that \(X\) is a positive semidefinite 
matrix. It serves as a simple and convex outer approximation of the \textit{Max-Cut polytope} 
\[
	\MC = \conv \{ X \in \SR ~|~ \rk X = 1\},
\]
where \(\conv\) denotes the convex envelope and \(\rk X\) denotes the rank of \(X\).

Note that \(\{ X \in \SR ~|~ \text{rk } X = 1 \} = \{ X ~|~ \exists x \in \{-1,1\}^n,~ X = x x^T\}\). Indeed a positive 
semidefinite matrix \(X\) has rank 1 if and only if there exists a nonzero vector \(x\) such that \(X = x x^T\). Then 
the condition \(\diag X = 1\)  implies that \(x_i^2 = 1\) for every \(i \in \{1, ..., n\}\), i.e., \(x_i = \pm 1\), and 
conversely.

The maximal value of a linear functional \( \langle A, . \rangle \) over a set \(E\) does not change if the set \(E\) is 
replaced by its convex envelope \(\conv E\). Therefore
\[
	\max_{
		\scriptsize
		\setlength{\tabcolsep}{0pt}
		\begin{array}{c}
			X = x x^T \\ 
			x \in \{-1,1\}^n
		\end{array}
	} \hspace{-3ex} \langle A, X \rangle = \max_{X \in \MC} \langle A, X \rangle.
\]

However, the Max-Cut polytope is a difficult polytope. Indeed, it has an exponential number of vertices and is defined 
by even more linear constraints. A good review of results on the Max-Cut polytope can be found in \cite{DezaLaurent}.

Maximizing \(\langle A, X \rangle \) over \(\SR\) instead of \(\MC\) for \(A \succeq 0\) approximates the exact solution 
of the problem with relative accuracy \(\mu = \dfrac{\pi}{2} - 1\) \cite{Nesterov}:
\[
	\dfrac{2}{\pi} \max_{X \in \SR} \langle A, X \rangle \leqslant \max_{X \in \MC} \langle A, X \rangle  
	\leqslant \max_{X \in \SR} \langle A, X \rangle.
\]

Define a function \(f: [-1, 1] \rightarrow [-1, 1]\) by \(f(x) = \dfrac{2}{\pi} \arcsin x\). Let \(\mathbf{f}\) be the 
operator which applies \(f\) element-wise to a matrix. A non-convex inner approximation of \(\MC\) is given by the 
\textit{trigonometric approximation}
\[
	\TA = \{ \mathbf{f}(X) ~|~ X \in \SR \}.
\]

Nesterov proved in \cite[Theorem 2.5]{Nesterov} that
\[
	\max_{X \in \TA} \langle A, X \rangle = \max_{X \in \MC} \langle A, X \rangle.
\]

Although not convex, \(\TA\) is simpler than \(\MC\) in the sense that checking whether a matrix \(X\) is in \(\TA\) 
can be done in polynomial time by computing \(\mathbf{f}^{-1}(X)\) and checking whether \(\mathbf{f}^{-1}(X)\) is 
in \(\SR\). This allows to reformulate the initial difficult problem (\ref{eqn:problem}) as an optimization problem over 
the algorithmically accessible set \(\TA\). The complexity of the problem in this form arises solely from the non-convexity 
of this set. 

Hirschfeld studied \(\TA\) in \cite[Section 4]{Hirschfeld}. In this work, we prove that \(\TA\) possesses an 
additional beneficial property. Namely, we prove the conjecture of Hirschfeld that it is starlike, i.e., for every 
\(X \in \TA\) and every \(\lambda \in [0, 1]\), the convex combination \(\lambda X + (1 - \lambda) I\) of \(X\) and 
the central point \(I\), the identity matrix, is in \(\TA\).

\section{Hirschfeld's conjecture}

\renewcommand{\arraystretch}{1.4}

In this section, we describe the conjecture and related results which have been obtained by Hirschfeld in his thesis 
\cite[Section 4.3]{Hirschfeld}.

In order to show that \(\TA\) is star-like, one has to prove that 
\[
	\forall X \in \SR, ~\mathbf{f}^{-1}(\lambda \mathbf{f} X + (1 - \lambda) I) \in \SR.
\]

\begin{sloppypar}
Note that the operator acting on X is nearly an element-wise one, defined by the function
\[
\begin{array}{c c c c}
	f_\lambda : & [-1, 1] & \longrightarrow & [-1, 1] \\
	& x & \longmapsto & f^{-1}(\lambda f(x)) = \sin(\lambda \arcsin x)
\end{array}
\]
acting on the off-diagonal elements, while the diagonal elements remain equal to 1, contrary to 
\({f_\lambda(1) = f^{-1}(\lambda) = \sin \dfrac{\pi \lambda}{2}}\). Thus one has to show that
\[
	\forall x \in \SR,~ \mathbf{f}_\lambda(X) + \left( 1 - \sin \dfrac{\pi \lambda}{2} \right) I \succeq 0.
\]
\end{sloppypar}

A sufficient condition is that \(\mathbf{f}_\lambda(X) \succeq 0\) for all \(X \in \SR\), i.e., the element-wise operator 
\(\mathbf{f}_\lambda\) is positivity preserving. Hirschfeld conjectured that this sufficient condition is verified
\cite[Conjecture 4.9]{Hirschfeld}.

\begin{lemma}
\label{lemma:pospres}
\[
	\forall X \in \SR,~ \mathbf{f}_\lambda(X) \succeq 0
\]
\end{lemma}

A sufficient (and necessary) condition for an operator of this type to be positivity preserving is that all of the Taylor 
coefficients of \(f_\lambda\) are nonnegative \cite{Schoenberg}.

Lemma \ref{lemma:pospres} proves the following theorem.

\begin{theorem}
	\(\TA\) is star-like.
\end{theorem}

\section{Proof of the conjecture}

In this section, we prove Lemma \ref{lemma:pospres}.

\begin{proof}[\textsc{Proof.}]
Let \(\lambda \in [0, 1]\) and write \(f_\lambda\) as a power series
\[
	f_\lambda(x) = \displaystyle \sum_{n \in \mathbb{N}}{a_n(\lambda) x^n}.
\]
The first two derivatives of \(f_\lambda\) are given by
\[
	f_\lambda'(x) = \dfrac{\lambda}{\sqrt{1 - x^2}} \cos(\lambda \arcsin x)
\]
and 
\[
	f_\lambda''(x) = \dfrac{x}{1 - x^2} \dfrac{\lambda \cos(\lambda \arcsin x)}{\sqrt{1 - x^2}} 
	- \dfrac{\lambda^2}{1 - x^2} \sin(\lambda \arcsin x).
\]
Hence \(f_\lambda\) is a solution on \((-1, 1)\) of the differential equation
\[
	(1 - x^2) f_\lambda'' - x f_\lambda' + \lambda^2 f_\lambda = 0.
\]
Therefore, the Taylor coefficients of \(f_\lambda\) verify the recurrence relation
\[
	(n + 2) (n + 1) a_{n + 2}(\lambda) - n (n - 1) a_n(\lambda) - n a_n(\lambda) + \lambda^2 a_n(\lambda) = 0
\]
which can be be re-expressed as
\begin{equation}
	a_{n + 2}(\lambda) = \dfrac{n^2 - \lambda^2}{(n + 2) (n + 1)} a_n(\lambda)
	\label{eqn:recurrence}
\end{equation}
with initial conditions
\[
	\left\{
	\begin{array}{c c c}
		a_0(\lambda) & = & 0 \\
		a_1(\lambda) & = & \lambda
	\end{array}
	\right..
\]
Given that \(\lambda \in [0, 1]\), a trivial induction shows that
\[
	\forall n \in \mathbb{N},~ a_n(\lambda) \geqslant 0.
\]
\end{proof}

Recursion (\ref{eqn:recurrence}) also proves that the roots of the polynomials \(a_n(\lambda)\) are located at 
\(0, \pm 1, ..., \pm n\) and are given by the polynomials \(\widetilde{P}_n(\lambda)\) \cite[eq. 4.23]{Hirschfeld}, 
as also conjectured by Hirschfeld.

\bibliographystyle{plain}

\end{document}